\documentclass[11pt]{amsart}
\usepackage{amsmath, amscd}
 
\theoremstyle{plain}
\newtheorem{thm}{Theorem}[section]

\newtheorem{pro}[thm]{Proposition}

\theoremstyle{definition}

\def \R {\mathbf{R}}
\def \Z {\mathbf{Z}}

\def \CR {\mathcal R}

\def \re {representation}
\def \a {\alpha}
\def \b {\beta}

\def \g {\gamma}

\def \o {\omega}
\def \ov {\overline}

\def \p {\phi}
\def \r {\rho}
\def \s {\sigma}

\def \th {\theta}

\def \O {\Omega}

\def \x {\times}

\begin{document}
 
\baselineskip.525cm
 
\title[Figure eight knot]
{The symplectic Floer homology of the figure eight knot}
%\thanks{Partially supported by NSF Grant DMS 9626166}
\author[Weiping Li]{Weiping Li}
\address{Department of Mathematics, Oklahoma State University \newline
\hspace*{.175in}Stillwater, Oklahoma 74078-0613}
\email{wli@math.okstate.edu}
\date{\today}

\begin{abstract}
In this paper, we compute the symplectic Floer homology of the 
figure eight knot. This provides first nontrivial knot with
trivial symplectic Floer homology. 
\end{abstract}
\maketitle

\section{Introduction\label{Intro}}

In \cite{li}, we generalized the Casson-Lin invariant \cite{lin} to the 
symplectic theory point of view. Our symplectic Floer homology of knots serves 
a new invariant for knots, and its Euler characteristic is half of the
signature of knots. 

We showed that the symplectic Floer homology
of the unknotted knot is trivial in \cite{li}. 
The natural question arises as whether
there is a nontrivial knot with trivial symplectic Floer homology.
We answer this question in this paper by computing the
symplectic Floer homology of the figure eight knot.

Although we know that the signature of the figure eight knot
$4_1 = \ov{\s_1 \s_2^{-1} \s_1 \s_2^{-1}}$ is zero, the signature does not
suffice to give the information of our finer invariant - 
the symplectic Floer homology.
For the square knot, we computed in \cite{li1} that the symplectic Floer
homology is nontrivial even though its signature is zero.
Our main result is the following.

\noindent{\bf Theorem} {\em
The symplectic Floer homology of the figure eight knot 
$4_1 = \ov{\s_1 \s_2^{-1} \s_1 \s_2^{-1}}$ is
\[ HF_i^{\text{sym}}(4_1) =
CF_i^{\text{sym}}(4_1) = 0, \ \ \ 
\text{for all $i \in {\Z}_4$}.\]}

To our knowledge, this is the first trivial symplectic
Floer homology involving nontrivial information. 
It is still an open question about if there is a non-homotopy
3-sphere with trivial instanton Floer homology. We wish to build the
relation between our symplectic Floer homology of knots \cite{li} and the
instanton Floer homology of homology 3-spheres \cite{fl2} through the
Dehn surgery technique. Using the calculation of the figure
eight knot, we hope to find an example of non-homotopy
3-sphere with trivial instanton Floer homology.

\section{The symplectic Floer homology}

\subsection{The symplectic Floer homology of braids}

We briefly recall our definition of the Floer homology of braids in this
subsection. See \cite{li} for more details.

For any knot $K = \ov{\b}$ with $\b \in B_n$, the braid group, 
the space ${\CR}(S^2 \setminus K)^{[i]}$ can be identified with the space
of $2n$ matrices $X_1 \cdots, X_n$, $Y_1, \cdots, Y_n$ in $SU(2)$ satisfying
\begin{equation} \label{zero}
\mbox{tr} (X_i) = \mbox{tr} (Y_i) = 0, \ \ \ \ \mbox{for $i =1, \cdots, n$},
\end{equation}
\begin{equation} \label{product}
X_1 \cdot X_2 \cdots X_n = Y_1 \cdot Y_2 \cdots Y_n .
\end{equation}
Note that $\pi_1(S^2 \setminus K)$ is generated by
$m_{x_i}, m_{y_i} (i =1, 2, \cdots, n)$
with one relation $\prod^n_{i=1} m_{x_i} = \prod^n_{i=1} m_{y_i}$.
There is a unique reducible
conjugacy class of {\re}s $s_K: \pi_1(S^3 \setminus K) \to U(1)$ 
such that
\[ s_K([m_{x_i}])  = s_K([m_{y_i}]) = \left[ \begin{array}{cc}
i & 0 \\
0 & - i  \end{array} \right]  .\]

Let ${\CR}^*(S^2 \setminus K)^{[i]}$ be the subset of
${\CR}(S^2 \setminus K)^{[i]}$ consisting of irreducible {\re}s.
Then ${\CR}^*(S^2 \setminus K)^{[i]}$ is a monotone symplectic
manifold of dimension $4n-6$ by Lemma 2.3 in \cite{li}. 
The symplectic manifold $(M, \o)$ is called {\em monotone} if
$\pi_2(M) = 0$ or if there exists
a nonnegative $\a \geq 0$ such that
$I_{\o } = \a I_{c_1}$ on $\pi_2(M)$, where  
$I_{\o }(u) = \int_{S^2} u^*(\o ) \in {\R}$ and
$I_{c_1}(u) = \int_{S^2} u^*(c_1) \in \Z$ for $u \in \pi_2(M)$.
The braid $\b$ induces a diffeomorphism 
$\p_{\b }:  {\CR}^*(S^2 \setminus K)^{[i]} \to
{\CR}^*(S^2 \setminus K)^{[i]}$. The induced diffeomorphism $\p_{\b }$ is 
symplectic, and the fixed point
set of $\p_{\b }$ is ${\CR}^*(S^3 \setminus K)^{[i]}$ 
(see Lemma 2.4 in \cite{li}).

Let $H: {\CR}^*(S^2 \setminus K)^{[i]} \x {\R} \to {\R}$
be a $C^{\infty}$ time-dependent Hamiltonian function with 
$H(x, s) = H(\p_{\b}(x), s+1)$. Let $X_s$ be the corresponding vector field
from $\o (X_s, \cdot) = dH_s( \cdot, s)$, and $\psi_s$ be the 
corresponding flow
\[\frac{d \psi_s}{ds} = X_s \circ \psi_s, \ \ \
\psi_0 = id.\]
Then we have $\psi_{s+1} \circ \p_{\b}^H = \p_{\b} \circ \psi_s$, where
$\p_{\b}^H = \psi_1^{-1} \circ \p_{\b}$.
Let $\O_{\p_{\b}}$ be the space of smooth paths $\a$ in
${\CR}^*(S^2 \setminus K)^{[i]}$ such that $\a (s+1) = \p_{\b}(\a (s))$. 
The symplectic action 
$a_H: \O_{\p_{\b}} \to {\R}/ 2 \a N{\Z}$ is given by
\[ da_H(\g) \xi = \int_0^1 \o (\dot{\g} - X_s(\g), \xi) ds. \]
So the critical points of $a_H$ are the fixed points of $\p_{\b}^H$.
For $x \in \text{Fix}(\p_{\b}^H)$, define $\mu (x) = \mu_u (x, s) \pmod {2N}$,
where $\mu_u$ is the Maslov index and $N=N(K)$ is the minimal value of 
the first Chern number of 
the tangent bundle of ${\CR}^*(S^2 \setminus K)^{[i]}$.
The integer $N(K)$ is a knot invariant.

Thus we have a
${\Z}_{2N}$-graded symplectic Floer chain complex:
\[CF_i^{\text{sym}} = 
\{ x \in \mbox{Fix}({\p}_{\b}) \cap {\CR}^*(S^2 \setminus K)^{[i]} :
\mu (x) = i\}, \ \ \ i \in {\Z}_{2N}. \] 
The following is Proposition 4.1 and Theorem 4.2 of \cite{li}.
\begin{thm} \label{main}
For a knot $K= \ov{\b}$ with the property that 
$\pi_2({\CR}^*(S^2 \setminus K)^{[i]}) = 0$ or $\a N(K) = 0$,
there is a well-defined ${\Z}$-graded 
symplectic Floer homology $HF_*^{\mbox{sym}}(\p_{\b})$.
The symplectic Floer homology 
$\{HF_i^{\mbox{sym}}(\p_{\b})\}_{i \in {\Z}_{2N}}$ is a knot invariant
and its Euler number is half of the signature of the knot (see \cite{li}).
\end{thm}

\subsection{The symplectic Floer homology of the figure eight knots}

The figure eight knot $4_1$ has the braid representative
$\s_1 \s_2^{-1} \s_1 \s_2^{-1}$. The knot $4_1$ has signature zero
since $4_1$ is equivalent (by an orientation preserving homeomorphism)
to its mirror image $\ov{4_1}$. So the figure eight knot is
amphicheiral. Also it is well-known that
the figure eight knot is not a slice knot, and represents
an element of order 2 in the knot cobordism group (see \cite{ro}).

We calculate the 
symplectic Floer homology of the figure eight knot
by identifying 
the fixed points
of the induced symplectic diffeomorphism in \S 2.1.

Let ${\CR}^*(S^2 \setminus 4_1)^{[i]}$ be the subset of
${\CR}(S^2 \setminus 4_1)^{[i]}$ consisting of irreducible representations.
Then ${\CR}^*(S^2 \setminus 4_1)^{[i]}$ can be also identified with
$(H_3 \setminus S_3)/SU(2)$ in Lin's notation \cite{lin}, i.e.,
the set of 6-tuple $(X_1, X_2, X_3, Y_1, Y_2, Y_3) \in SU(2)^6$
satisfying
$\text{tr} (X_j) = \text{tr} (Y_j) = 0 (j = 1, 2, 3)$ and
\[ X_1 X_2 X_3 = Y_1 Y_2 Y_3.\]
By operating the conjugation on $X_3$ and $Y_3$, we may assume that
\[ X_3 = \begin{pmatrix}
i & 0\\
0 & -i \\
\end{pmatrix}, \ \ \
Y_3 = \begin{pmatrix}
i \cos \th  & \sin \th \\
- \sin \th  & - i \cos \th \\
\end{pmatrix}, \ \ \ 0 \leq \th \leq \pi .\]
If $\th = 0$ and $\pi$, then we get two copies of
$(H_2 \setminus S_2)/SU(2)$ which is the pillow case (a 2-sphere with
four cone points deleted \cite{li, lin}). 
For $0 < \th < \pi$, the identification reduces down to the following  
\[X_1 X_2
\begin{pmatrix}
\cos \th & - i \sin \th \\
- i \sin \th & \cos \th \\
\end{pmatrix} = Y_1 Y_2. \]
Let $R_{\th}$ be the representations in ${\CR}^*(S^2 \setminus 4_1)^{[i]}$
satisfying the above equation. So the space $R_{\th }$ is the non-singular
piece in ${\CR}^*(S^2 \setminus K)^{[i]}$. For $0 < \th , \th^{'} < \pi$,
the space $R_{\th }$ is diffeomorphic to the space $R_{\th^{'}}$.
In particular, they are all diffeomorphic to $R_{\pi/2}$.
In this case, we see that ${\CR}^*(S^2 \setminus 4_1)^{[i]}$ is a generalized
pillow case:
\[{\CR}^*(S^2 \setminus 4_1)^{[i]} = \bigcup_{0 \leq \th \leq \pi} R_{\th} .\]

The fixed point set of $\p_{4_1}$ is ${\CR}^*(S^3 \setminus 4_1)^{[i]}$
by Lemma 2.4 in \cite{li}. So we have, for $\s = \s_1 \s_2^{-1} 
\s_1 \s_2^{-1}$,
$\text{Fix}\, (\p_{4_1}) = \{(X_1, X_2, X_3) \in SU(2)^3|
\s (X_j) = X_j, j=1, 2, 3\}$ up to conjugation.
Let $B_n$ be the braid group of rank $n$ with the standard generators
$\s_1, \cdots, \s_{n-1}$, and $F_n$ be the free group of rank $n$
generated by $x_1, \cdots, x_n$. Then the automorphism of $F_n$ representing
$\s_k$ is given by (still denote it by $\s_k$)
\begin{equation} \label{action}
\begin{aligned}
\s_k : \ \ \ & x_k \mapsto x_k x_{k+1} x_k^{-1} \\
       & x_{k+1} \mapsto x_k \\
       & x_l \mapsto x_l, \ \ \ l \neq k, k+1.
\end{aligned} 
\end{equation}

By (\ref{action}), we compute the followings.
\[\begin{aligned}
\s_1 \s_2^{-1} \s_1 \s_2^{-1} (x_1) & = 
\s_1 \s_2^{-1} \s_1 (x_1^{-1}) = \s_1 \s_2^{-1} (x_1x_2^{-1}x_1^{-1}) \\
& = \s_1 (x_1 x_2 x_3 x_2^{-1} x_1^{-1}) \\
& = (x_1 x_2 x_1^{-1})x_1 x_3 x_1^{-1} (x_1 x_2 x_1^{-1})^{-1}\\
& = x_1 x_2 x_3 x_2^{-1} x_1^{-1}. \\
\s_1 \s_2^{-1} \s_1 \s_2^{-1} (x_2) & = 
\s_1 \s_2^{-1} \s_1 (x_2x_3^{-1}x_2^{-1}) \\
& = \s_1 \s_2^{-1}(x_1x_3^{-1}x_1^{-1}) = \s_1(x_1^{-1}x_2x_1) \\
& = x_1 x_2^{-1} x_1 x_2 x_1^{-1}.\\
\s_1 \s_2^{-1} \s_1 \s_2^{-1} (x_3) & = 
\s_1 \s_2^{-1} \s_1 (x_2^{-1}) = \s_1 \s_2^{-1}(x_1^{-1})\\
& = \s_1(x_1) = x_1 x_2 x_1^{-1}.
\end{aligned}\]

Therefore the fixed point set of $\p_{4_1}$ is the set of points
$(X_1, X_2, X_3) \in SU(2)^3$ such that
\[\begin{aligned}
 \text{tr} (X_j) & = 0, \ \ \ j=1, 2, 3,\\
X_1 X_2 X_3 X_2^{-1} X_1^{-1} & = X_1,\\
X_1X_2^{-1} X_1 X_2 X_1^{-1}& = X_2, \\
 X_1 X_2 X_1^{-1} & = X_3, \end{aligned} \]
up to conjugation. Up to conjugation, we can assume that
\[X_2 = \begin{pmatrix}
i & 0\\
0 & -i \\
\end{pmatrix}, \ \ \
X_1 = \begin{pmatrix}
i \cos \th  & \sin \th \\
- \sin \th  & - i \cos \th \\
\end{pmatrix}, \ \ \ 0 \leq \th \leq \pi .\]
From the last equation in the above, we obtain
\[\begin{aligned}
X_1 X_2 X_1^{-1} & = 
\begin{pmatrix}
i \cos \th  & \sin \th \\
- \sin \th  & - i \cos \th \\
\end{pmatrix} \begin{pmatrix}
i & 0\\
0 & -i \\
\end{pmatrix} \begin{pmatrix}
- i \cos \th  & - \sin \th \\
\sin \th  & i \cos \th \\
\end{pmatrix}  \\
& = \begin{pmatrix}
i \cos 2 \th  & \sin 2 \th \\
- \sin 2 \th  & - i \cos 2 \th \\
\end{pmatrix} = X_3.
\end{aligned} \]
So the matrix $X_3$ is completely determined by the parameter
$\th \in [0, \pi]$. This is, in fact, a key to complete the calculation.
Now substituting $X_3$ into the relation $\s_1 \s_2^{-1} \s_1 \s_2^{-1}(X_1)
= X_1$, we have
\[\begin{aligned}
X_1 X_2 X_3^{-1} X_2^{-1} X_1^{-1} & = 
\begin{pmatrix}
- \cos \th  & - i \sin \th \\
- i \sin \th & - \cos \th \\ \end{pmatrix}
\begin{pmatrix}
- i \cos 2 \th  & \sin 2 \th \\
- \sin 2 \th & i \cos 2 \th \\
\end{pmatrix} 
\begin{pmatrix} 
- \cos \th  & i \sin \th \\
i \sin \th & - \cos \th \\ \end{pmatrix} \\
& = \begin{pmatrix}
- i \cos 4 \th  & - \sin 4 \th \\
\sin 4 \th & i \cos 4 \th  \\
\end{pmatrix} 
= X_1. \end{aligned} \]
This reduces to the equations
\begin{equation} \label{4th}
\cos 4 \th = - \cos \th , \ \ \ \ \ 
\sin 4 \th = - \sin \th. 
\end{equation}
Similarly, we compute
\[\s_1 \s_2^{-1} \s_1 \s_2^{-1} (X_2)  = 
\begin{pmatrix}
i \cos 3 \th  &  \sin 3 \th \\
- \sin 3 \th & - i \cos 3 \th  \\ \end{pmatrix} 
= X_2 = \begin{pmatrix}
i & 0\\
0 & -i \\
\end{pmatrix}, \]
to get the equations
\begin{equation} \label{3th}
\cos 3 \th = 1, \ \ \ \ \ \sin 3 \th = 0.
\end{equation}
Thus the fixed point of $\p_{4_1}$ can be identified with
\[X_1 = \begin{pmatrix}
i \cos \th  & \sin \th \\
- \sin \th  & - i \cos \th \\
\end{pmatrix}, 
X_2 = \begin{pmatrix}
i & 0\\
0 & -i \\
\end{pmatrix}, 
X_3 = \begin{pmatrix}
i \cos 2 \th  & \sin 2 \th \\
- \sin 2 \th  & - i \cos 2 \th \\
\end{pmatrix},  \ \ 0 \leq \th \leq \pi ,\]
subject to equations (\ref{4th}) and (\ref{3th}). Using the equations
(\ref{3th}) and the angle addition formulae for sine and cosine
functions with $4 \th = 3 \th + \th$, (\ref{4th}) becomes
\begin{equation} \label{th}
\sin \th = 0, \ \ \ \cos \th = 0.
\end{equation}
There is no solution for (\ref{th}). Hence
\begin{equation} \label{epty}
\text{Fix}\, (\p_{4_1}) = \emptyset \ \ \ \text{(empty set)}.
\end{equation}
 
\begin{thm} 
The symplectic Floer homology of the figure eight knot
$4_1 = \ov{\s_1 \s_2^{-1} \s_1 \s_2^{-1}}$ is
\[ HF_i^{\text{sym}}(4_1) =
CF_i^{\text{sym}}(4_1) = 0, \ \ \
\text{for all $i \in {\Z}_{2N}$}.\]
\end{thm}
Proof: Since the $\Z_{2N}$-graded symplectic Floer
chain complex $CF_i^{\text{sym}}(4_1)$ is generated
by $\text{Fix}\, (\p_{4_1})$, the result follows from (\ref{epty}).
\qed

\subsection{The symplectic Floer homology of knots with braid representatives
in $B_3$}

It seems that the method in \S 2.2 can be adapted to knots with 
braid representatives in $B_3$. We are going to illustrate another example
to show that the computation for the figure eight knot in \S 2.2 is quite
lucky.

Let $K = 5_2$ be the knot with 5-crossings. We have the braid representative
$\s_1^2 \s_2^2 \s_1^{-1} \s_2$ for the knot $5_2$ (see \cite{ro}). Thus the
fixed points of $\p_{5_2}$ can be identified, by the same method in \S 2.2,
with the set of points $(X_1, X_2, X_3) \in SU(2)^3$ such that
\[\begin{aligned}
 \text{tr} (X_j) & = 0, \ \ \ j=1, 2, 3,\\
X_1 X_2 X_3 X_1 X_2^{-1} X_1^{-1}X_3^{-1} X_2^{-1} X_1^{-1} & = X_1,\\
X_1X_2 X_3^{-1} X_1^2 X_2^{-1} X_1^{-1}& = X_2, \\
 X_1 X_2 X_1^{-1} X_2^{-1} X_1^{-1} & = X_3, \end{aligned} \]
up to conjugation. This follows a straightforward calculation of
$\s_1^2 \s_2^2 \s_1^{-1} \s_2 (x_j) (j=1, 2, 3)$. Again we can compute $X_3$
from the last equation in the above.
\[X_1 X_2 X_1^{-1} X_2^{-1} X_1^{-1}  =
\begin{pmatrix}
- i \cos 3 \th  & - \sin 3 \th \\
 \sin 3 \th  &  i \cos 3 \th \\
\end{pmatrix} 
= X_3. \]
Then $\s_1^2 \s_2^2 \s_1^{-1} \s_2 (X_j) = X_j (j=1, 2)$ gives us
\[\begin{aligned}
\s_1^2 \s_2^2 \s_1^{-1} \s_2 (X_1) & =
\begin{pmatrix}
- i \cos 6 \th  & - \sin 6 \th \\
 \sin 6 \th  &  i \cos 6\th \\
\end{pmatrix} = X_1 \\
\s_1^2 \s_2^2 \s_1^{-1} \s_2 (X_2) & =
\begin{pmatrix}
- i \cos 5 \th  & - \sin 5 \th \\
 \sin 5 \th  &  i \cos 5 \th \\
\end{pmatrix} = X_2.
\end{aligned} \]
Thus we need to solve the equations
\begin{equation} \label{star}
\cos 6 \th  = - \cos \th, \ \ \ \    \sin 6 \th = - \sin \th, 
\cos 5\th  =  - 1, \ \ \ \  \sin 5 \th = 0.
\end{equation}
There are three solutions of (\ref{star}) with $\th =  \frac{\pi}{5},
\frac{3\pi}{5}, \pi$. Let $\r_j (j=1, 2, 3)$ be the corresponding fixed
points of $\p_{5_2}$ in ${\CR}^*(S^2 \setminus 5_2)^{[i]}$.

By following the method in \cite{gl}, for $K = 5_2$, we have all
type I double points so that the correction term $\mu = 0$. Using the
definition of Goeritz matrix in \S 1 of \cite{gl}, we get the
Goeritz matrix of $5_2$:
\[ G(5_2) = \begin{pmatrix}
4 & -3 & -1 \\
-3 & 4 & -1 \\
-1 & -1& 2 \\ \end{pmatrix}.\]
By the theorem 6 of \cite{gl}, we have
\[\text{Signature}(5_2) = \text{Signature}(G(5_2)) - \mu = 2.\]
By Theorem~\ref{main}, the Euler characteristic of the symplectic Floer
homology of $5_2$ is one.
\begin{pro}
The symplectic Floer chain complex of $5_2$ is given by:
one of the odd chain groups is generated by one of $\r_j (j=1, 2, 3)$;
even chain groups are generated by the rest two fixed points of
$\p_{5_2}$.
\end{pro}

It is nontrivial to determine the Maslov index of $\r_j$ and the
possible Floer boundary map in order to complete the calculation.

%\noindent{\bf Acknowledgment:} 

\end{document}